 \documentclass[11pt,letterpaper]{amsart}
\pdfoutput=1
\usepackage[T1]{fontenc}
\usepackage{lmodern}

\usepackage{amsmath,amsfonts,amssymb,amsthm,stmaryrd,bbm,graphicx,mathtools,enumerate}
\usepackage{mathrsfs}
\usepackage[T1]{fontenc} 
\usepackage[latin1]{inputenc}
\usepackage[english]{babel}
\usepackage[tracking,spacing,kerning,babel]{microtype}
\usepackage{wasysym} 
\usepackage{color}
\usepackage{xcolor}
\usepackage{framed} 

\usepackage{wrapfig} 

\usepackage[final]{hyperref}   
\hypersetup{
    linktoc=page,
    linkcolor=red,          
    citecolor=blue,        
    filecolor=blue,      
    urlcolor=cyan,
   colorlinks=true           
}

\usepackage{lipsum} 

\usepackage[toc,page]{appendix}

\usepackage{minitoc}
\usepackage{multicol}

\newtheoremstyle{mine}
{\baselineskip}
{\baselineskip}
{\itshape}
{
}
{\bfseries}
{.}
{.5em}
{#1 #2\ifx#3\relax\else~(#3)\fi}

\theoremstyle{mine}

\newtheorem{theorem}{Theorem}
\numberwithin{theorem}{section}
\newtheorem{corollary}[theorem]{Corollary}
\newtheorem{proposition}[theorem]{Proposition}
\newtheorem{lemma}[theorem]{Lemma}

\newtheorem{definition}[theorem]{Definition}

\numberwithin{equation}{section}

\theoremstyle{remark}
\newtheorem{remark}{Remark}

\colorlet{shadecolor}{blue!10}

\def\rm{\reversemarginpar}

\let\qed=\QED

\renewcommand{\epsilon}{\varepsilon}

\newcommand{\R}{\mathbb{R}}

\newcommand{\Z}{\mathbb{Z}}
\newcommand{\N}{\mathbb{N}}
\renewcommand{\S}{\mathbb{S}}

\def\T{\mathbb{T}}

\def\calE{\mathcal{E}}
\def\calF{\mathcal{F}}

\def\Cov#1{\mathrm{Cov}\bigl[ #1\bigr]}

\def\Im{{\rm Im}\,}
\def\Re{{\rm Re}\,}

\def\P{\mathbb{P}} 
\def\E{\mathbb{E}} 
\def\md{\mid}

\def \eps {\epsilon}

\def\Bb#1#2{{\def\md{\bigm| }#1\bigl[#2\bigr]}}

\def\Eb{\Bb\E}

\def\<#1{\langle #1\rangle}

\definecolor{darkgreen}{rgb}{0,0.6,0.05}

\def\nn{\nonumber}
\def\bi{\begin{itemize}}  
\def\ei{\end{itemize}}
\def\bnum{\begin{enumerate}} 
\def\enum{\end{enumerate}}
\def\ni{\noindent}
\def\bf{\bfseries}

\author{Christophe Garban}
\author{Vincent Vargas}

\title[Harmonic analysis of Gaussian multiplicative chaos]{Harmonic analysis of Gaussian multiplicative chaos on the circle}

\address
{Universit\'e Claude Bernard Lyon 1, CNRS UMR 5208, Institut Camille Jordan, 69622 Villeurbanne, France \,, Institut Universitaire de France (IUF)}
\email{garban@math.univ-lyon1.fr}

\address{Universit\'e de Gen\`eve, Section de math\'ematiques, UNI DUFOUR, 24 rue du G\'en\'eral Dufour,
CP 64 1211 Geneva 4, Switzerland}
\email{Vincent.Vargas@unige.ch}

\begin{document}

\maketitle

\begin{abstract}
In this paper, we initiate the harmonic analysis of Gaussian multiplicative chaos (GMC) on the circle, i.e. the study of its Fourier coefficients. In particular, we show that almost surely GMC is a so-called Rajchman measure which means that its Fourier coefficients converge to $0$ when the frequency goes to infinity. We supplement this result with a convergence in law result for the rescaled Fourier coefficients.
\end{abstract}

\section{Introduction}
\subsection{Introduction and state of the art.}

Given a Radon measure $\nu(d \theta)$ on the standard unit cercle $[0,2\pi)$, it can be a notoriously difficult problem to study the asymptotics of $c_n:= \int_0^{2\pi}  e^{i n \theta} \nu(d\theta)$ as $n$ goes to infinity, especially when $\nu$ lives on a fractal set $F$ (i.e. a set of Hausdorff dimension strictly less than $1$).
In this context, $\nu(d \theta)$ is a so-called  {\em Rajchman measure} if   $c_n$ goes to $0$ as $n$ goes to infinity. The importance of this property is intimitaly related to the concept of {\em set of uniqueness} whose definition goes back to Cantor during the 19th century\footnote{Cantor's study of sets of uniqueness highly inspired his introduction of set theory.} and later studied by Lebesgue among others: see  \cite{lyons2020seventy,korner1992sets,lyons1985fourier} for a very nice account on this and a precise definition of a set of uniqueness. Roughly speaking, a set of uniqueness is a set which is ``small'' from the point of view of Fourier analysis. Rajchman's property is intimately linked to sets of uniqueness as follows. If $F$ is a closed fractal set which may be equipped with a Rajchman measure $\nu$ supported on $F$, then $F$ is not a set of uniqueness  (in which case $F$ is called a set of multiplicity). The main purpose of this article is precisely to prove that a canonical multifractal (random) measure called Gaussian multiplicative chaos (GMC) is a.s. a Rajchman measure. We expect our methodology can be adapted to other natural multifractal measures on the circle like the ones studied in \cite{Fan_Meyer}.

Gaussian multiplicative chaos (GMC), introduced in the pioneering work of Kahane \cite{Kahane}, is a theory of random measures which has received a lot of attention recently in probability theory, due to its relation to finance, turbulence, conformal field theory and Schramm-Loewner evolution/conformal loop ensembles. In dimension $d$ and for a real parameter $\gamma \geq 0$, GMC measures are formally defined on an open set $O \subset \R^d$ as
\begin{equation}\label{measintrorev}
M_\gamma (dx)= e^{\gamma X(x)-\frac{\gamma^2}{2} \E[X(x)^2]} dx
\end{equation}
where $dx$ is the Lebesgue measure and $X(x)$ is a centered Gaussian field which is log-correlated
\begin{equation*}
\E[ X(x) X(y) ]= \log \frac{1}{|x-y|}+f(x,y)
\end{equation*}
where $f$ is some smooth function on $O \times O$. The covariance kernel thus possesses a singularity along the diagonal and it is clear that
making sense of \eqref{measintrorev} is not straightforward since $X$ is not defined pointwise but rather lives with probability $1$ in the space of distributions (in the sense of Schwartz). The standard approach consists in applying an ultraviolet regularization to the distribution $X$ in order to get rid of the singularity of the covariance kernel. The regularization usually depends on a small parameter, call it $\eps$, that stands for the extent to which the field has been regularized. The measure \eqref{measintrorev} is naturally understood as the 
limit of the random measures:
\begin{equation}\label{limitreal}
M_{\gamma,\eps}(A) =\int_Ae^{\gamma X_\eps(x)-\frac{\gamma^2}{2}\E[X_\eps^2(x)]} dx
\end{equation}
 when the regularization parameter $\eps$ goes to $0$. It is known since the work of Kahane \cite{Kahane} that this produces non trivial limiting objects when the real parameter $\gamma$ is strictly less than the critical value $\gamma_c=\sqrt{2d}$.  

%

The simplest and most symmetrical model of GMC is undoubtedly the GMC built from the Gaussian Free Field (GFF) $\varphi$ on the standard unit circle $[0,2\pi)$ with covariance kernel   
\begin{equation*}
\E[\varphi(\theta)\varphi(\theta')]= \ln \frac{1}{|e^{i\theta}-e^{i\theta'} |}.
\end{equation*}

From now on, we will consider this case of GMC in the sequel. We will work with the cut-off

\begin{equation*}
\varphi_\epsilon:= \frac{1}{2}\int_{-1}^1 \varphi (\theta+\epsilon v) dv
\end{equation*}
with covariance given by
\begin{equation*}
\E[\varphi_\epsilon(\theta)\varphi_\epsilon(\theta')]=  \frac{1}{4}\int_{-1}^1 \int_{-1}^1 \ln \frac{1}{|e^{i(\theta+\epsilon v)}-e^{i(\theta'+\epsilon v')} |} dv dv'
\end{equation*}
which converges to $\ln \frac{1}{|e^{i\theta}-e^{i\theta'} |}$ when $\epsilon$ goes to $0$ and such that $\E[\varphi_\epsilon(\theta)^2]= \ln \frac{1}{\epsilon}+C+o(1)$ for some constant $C$.

The GMC measure is defined as the limit in probability in the space of Radon measures (see \cite{NathanaelGMC} for an elegant approach)
\begin{equation}\label{defconvergencechaos}
M_\gamma(d \theta):=e^{\gamma \varphi (\theta)-\frac{\gamma^2}{2} \E[\varphi(\theta)^2 ]}d\theta := e^{-\frac{\gamma^2}{2} C} \underset{\epsilon \to 0}{\lim} \;  \epsilon^{\gamma^2/2} e^{\gamma \varphi_\epsilon(\theta) } d\theta
\end{equation}
 In the sequel, we suppose that $\gamma$ belongs to $[0, \sqrt{2})$ since for $\gamma \geq 2$ the limiting measure $M_\gamma(d \theta)$ is $0$. We will also set $M_{\gamma,\epsilon}(d\theta)=  \epsilon^{\gamma^2/2} e^{\gamma \varphi_\epsilon(\theta) } d\theta $. Since the GMC $M_\gamma$ lives on the circle, one can define its Fourier coefficients
\begin{equation}\label{Fourier}
c_n= c_n(M_\gamma)= \int_0^{2\pi}  e^{i n \theta} e^{\gamma \varphi (\theta)-\frac{\gamma^2}{2} \E[\varphi(\theta)^2 ]}d\theta
\end{equation}
In addition to the fact that it is very natural to study the coefficients $c_n$ within the framework of harmonic analysis, they are also an essential ingredient in the construction of the Virasoro algebra of Liouville conformal field theory \cite{BGKRV} where they appear as potential terms reflecting the interacting nature of the theory (as compared to Gaussian Free Field theory). Moreover, we expect that the sequence $c_n$ exhibits the same behaviour as a wide range of natural complex valued partial sums of random multiplicative products: see below for a discussion on the expected relation between the $c_n$ and the other settings which appear in number theory or random matrix theory. Also in relation with the study of random matrix theory and closely related to our framework is the work \cite{CN} where the authors study and derive the explicit law of the so-called {\em Verblunsky coefficients} of $M_\gamma$. These coefficients appear in the theory of orthogonal polynomials associated to $M_\gamma$ and are closely related to the Fourier coefficients $c_n$. However, it appears that this relation is not sufficiently explicit to deduce our results from their results.

%
%

It is well known that $M_\gamma$ is a multifractal measure, i.e. it is not scale invariant but rather possesses a continuum spectrum of scaling behaviors and lives on a set $E_\gamma$ of Hausdorff dimension $1-\frac{\gamma^2}{2}$: see \cite{RV} for instance for an account on this. In fact, $M_\gamma$ is exactly  of dimension $1-\frac{\gamma^2}{2}$ since $M_\gamma(d \theta)$ almost everywhere 
\begin{equation*}
\underset{r \to 0}{\lim} \: \frac{M_\gamma[\theta-r, \theta+r]}{\ln r}  = 1-\frac{\gamma^2}{2}
\end{equation*} 
One can choose $E_\gamma$ as the set of points $\theta$ with the above behaviour and this also implies that the Hausdorff dimension $\text{dim}_H(M_\gamma)$ of the measure is $1-\frac{\gamma^2}{2}$ (which means that all sets of dimension striclty less than $1-\frac{\gamma^2}{2}$ do not charge $M_\gamma$: see \cite{Falk}).

As mentioned in the beginning of the introduction, a very natural question is the behaviour of $c_n$ as $n$ goes to infinity and in particular if $M_\gamma$ is a Rajchman measure, i.e. does $c_n$ go to $0$ as $n$ goes to infinty? This has been a long standing question on the measure $M_\gamma$ and the main result of this paper is to settle this question in the affirmative.

\subsection{The main results and perspectives}

\subsubsection{Main results.}

The first main result of the paper is the fact that $M_\gamma$ is a Rajchman measure: 

\begin{theorem}\label{Rajchmanback}
For all $\gamma< \sqrt{2}$, the sequence $c_n$ goes to $0$ almost surely as $n$ goes to infinity.
\end{theorem}


\begin{remark}
One may also consider GMC on the standard $d$ torus $[0,2\pi)^{d}$ for any $d \geq 1$. In the case $d=2$ (of particular importance), one considers the exponential of the GFF whose covariance is given also by the Green function on the torus. In the general case, one considers the exponential of the field $X$ defined via the following series which converges a.s. in the space of distributions  
\begin{equation*}
X(\theta)= \sum_{n=(n_1, \cdots,n_d) \in \N^{d} \setminus \lbrace 0 \rbrace}  \frac{e^{i n \cdot \theta}  \epsilon_n}{|n|^{d/2}}, \quad \theta \in [0,2\pi)^{d}
\end{equation*}
where $|n|:= \sqrt{n_1^2+ \cdots n_d^2} $, $n \cdot \theta$ denotes the standard scalar product and $(\epsilon_n)$ is an i.i.d. sequence of standard complex Gaussian variables.
One can then define the Fourier coefficients $c_n$ as  $c_n=int_{[0,2\pi)^{d}} e^{i n \theta} e^{\gamma \varphi (\theta)-\frac{\gamma^2}{2} \E[\varphi(\theta)^2 ]}d\theta$ where $d\theta$ denotes the standard Lebesgue measure, which is defined via the theory of GMC (see ) and is non trivial for $\gamma^2<2d$. Our methods are robust in the sense that with minor modifications one can show that for $\gamma^2<2d$, the sequence $c_n$ goes to $0$ when $|n|$ goes to infinity. An analogue of Theorem \ref{RajchmanbackL4} and of Theorem \ref{thlaw} can also be proved under the condition $\gamma< \sqrt{\frac{d}{2}}$.
\end{remark}

To the best of our knowledge, the only instance where such a result is proved on a (random) multifractal measure is the paper \cite{FalkJin} where such a property is proved for a 2d version of GMC  for small values of $\gamma$ (roughly they consider the case $\gamma<0.25$). In this case one can even find a (non explicit) bound for the Fourier dimension of the 2d measure considered. Recall that the Fourier dimension of $M_\gamma$ is by definition the maximal $s$ such that $|c_n| = O(\frac{1}{n^{s/2}} )$ as $n$ goes to infinity. When $\gamma < \frac{1}{\sqrt{2}}$, one can reinforce Theorem \ref{Rajchmanback} by a simple method of moments and get (clearly non optimal) bounds on the Fourier dimension:

\begin{theorem}\label{RajchmanbackL4}
For all $\gamma< \frac{1}{\sqrt{2}}$, for all $\beta \in (0,\frac{1-2 \gamma^2}{4})$ there exists a random constant $C$ such that almost surely
\begin{equation*}
|c_n| \leq \frac{C}{n^{\beta}}
\end{equation*}
\end{theorem}
Therefore the Fourier dimension is bounded below in this case by $\frac{1-2 \gamma^2}{2}$. As mentioned above, we do not expect Theorem \ref{RajchmanbackL4} to be optimal: see the discussion below where we would rather expect that the Fourier dimension is $1-\gamma^2$  for $\gamma< \frac{1}{\sqrt{2}}$.

In view of Theorem \ref{Rajchmanback}, it is natural to seek for a renormalization of  $c_n$ such that convergence to a non trivial variable holds. This will be the object of the next theorem for the case $\gamma< \frac{1}{\sqrt{2}}$. First, for any $\gamma<1$, let us define 
\begin{equation*}
\kappa=\kappa(\gamma):= \int_\R \frac{e^{i 2\pi v}}{|v|^{\gamma^2}} dv
\end{equation*}

We have the following convergence in law:
\begin{theorem}\label{thlaw}
For $\gamma< \frac{1}{\sqrt{2}}$ we have
\begin{equation*}
n^{(1-\gamma^2)/2} c_n \underset{n \to \infty}{\overset{(d)} {\rightarrow}} \sqrt{\frac{\kappa}{2}} \,   W_{M_{2\gamma}[0,2\pi ] }
\end{equation*}
where  $M_{2\gamma}$ is the GMC of parameter twice $\gamma$ and $W$ is a complex Brownian motion independent of $M_{2\gamma}$.
\end{theorem}


\subsubsection{Perspectives and conjectures.}

The above results bring a few comments and natural open questions. Our first remark deals with the generality of the methodology behind the above results: 

\begin{remark}
The scaling in theorem Theorem \ref{thlaw} can be understood from the point of view of capacity estimates of $M_\gamma$. Recall the following result from \cite{LRV}
\begin{equation}\label{capacityestimate}
\int_{[0,2\pi] \times [0, 2 \pi]} \frac{M_\gamma (d \theta) M_\gamma (d \theta')}{|e^{i \theta}- e^{i \theta'}|^s} < \infty, \quad a.s.
\end{equation}
if and only if \footnote{The exponent $s$ used here corresponds to $\beta^2$ in \cite{LRV}.}
\begin{equation*}
\begin{cases}
s < s_\star(\gamma):= 1-\gamma^2 \: \text{if} \: \gamma < \frac{1}{\sqrt 2}. \\
s < s_\star(\gamma):=(\sqrt 2 -\gamma) ^2 \: \text{if} \: \gamma \geq \frac{1}{\sqrt 2}. \\
\end{cases}
\end{equation*}

A standard capacity computation, which can be obtained by developing $M_\gamma$ along the Fourier basis, yields the following bound for all $s \in (0,1)$
\begin{equation}\label{capacityequivalent}
c_s  \sum_{n \in \Z} n^{s-1}|c_n|^2 \leq \int_{[0,2\pi] \times [0, 2 \pi]} \frac{M_\gamma (d \theta) M_\gamma (d \theta')}{|e^{i \theta}- e^{i \theta'}|^s} \leq C_s  \sum_{n \in \Z} n^{s-1}|c_n|^2
\end{equation}
where $c_s,C_s>0$ are two deterministic constants which depend only on s (and not on $M_\gamma$). Clearly in view of the above discussion the Fourier dimension of $M_\gamma$ is smaller than $s_\star(\gamma)$ and strictly smaller than the Hausdorff dimension $\text{dim}_H(M_\gamma)$ since in this setting one can notice the strict inequality $s_\star(\gamma)< \text{dim}_H(M_\gamma)$ (the fact that $s_\star(\gamma) \leq \text{dim}_H(M_\gamma)$  is a general property of Radon measures). Hence it is very natural to ask: is the Fourier dimension of $M_\gamma$  given by $s_\star(\gamma)$ (in which case $M_\gamma$ is a so-called Salem measure)? 
\end{remark}


\begin{remark}
Finally, it is possible to define $M_\gamma$ for $\gamma=\gamma_c=\sqrt{2}$ by using a slightly different limiting procedure: this is the so-called derivative martingale which lives on a set of Hausdorff dimension $0$ but has no atoms. See for example \cite{powell2020critical, aru2019critical} and references therein. Our method is not applicable to this case and it would be very interesting to prove that $M_{\gamma_c}$ for $\gamma_c=\sqrt{2}$ is also almost surely a Rajchman measure, which would then have Fourier dimension $0$: see \cite{korner1992sets} for statements on such measures which satisfy the Rajchman property.   
\end{remark} 

In conclusion, to summarize the above discussion, here is a list of open problems which seem interesting:

\vspace{0.1 cm}

\noindent
{\bf Open problems.}

\begin{itemize}
\item
Show the optimal asymptotic of $c_n$ in the almost sure sense. (Theorem \ref{thlaw} only gives a strong hint of what this a.s. optimal asymptotic should be when $\gamma < \frac {1} {\sqrt{2}}$). 
\item
When are the convolutions of $M_\gamma$ H\"older functions?
\item
How can one prove Rajchman's property for the critical case $\gamma_c=\sqrt{2}$?
(This question looks particularly interesting to us as our present technique would somehow require to consider an infinite convolution $M_{\gamma_c}^{*\infty}$ which is of course ill-defined). 
\item Is it the case that the Fourier dimension of $M_\gamma$ is a.s. given by $s_\star(\gamma)$? (Already for $\gamma <\frac {1} {\sqrt{2}}$, this question is open). 
\end{itemize}

\subsubsection{Relation to other models and results.}

To the best of our knowledge, Theorem \eqref{Rajchmanback} is the first proof of the Rajchman property for a (random) multifractal measure except for the work \cite{FalkJin} which considers small values of $\gamma$ for a 2d GMC. Let us note that in \cite{ShmerSuo}, the authors prove such a property for a monofractal random measure which roughly corresponds to replacing the normal field $X$ by a random variable taking values in a set of cardinality two; in this case, they even determine the exact rate of decay. However, though the methods of \cite{ShmerSuo} could perhaps be adapted to handle certain multifractal measures, the setting of this paper seems out of reach of the techniques developed in \cite{ShmerSuo}. In particular, in this monofractal setting the Fourier dimension coincides with the Hausdorff dimension of the measures whereas by the above discussions we know that this not true for $M_\gamma$.

\vspace{0.1 cm}
We first discuss two related models which have already been mentioned above (namely {\em Verblunsky coefficients}  and  Liouville conformal field theory) and will discuss further below two other related models respectively in analytic number theory and random matrix theory.

\bi
\item In \cite{CN}, the authors analyse the so-called {\em circular $\beta$ ensemble} as the zeros of random orthogonal polynomials specified by an explicit law on their {\em Verblunsky coefficients}. The random measure corresponding to these Verblunsky coefficients corresponds to a GMC like measure on the circle. The relationship with our work reads as follows: if $\mu$ is a measure on the circle $\T^1$ which is a small perturbation of the Lebesgue measure ($\mu(dx) = dx + \eps \nu$), then its $j^{th}$ Verblunsky coefficient $\alpha_j$ is related to the $j^{th}$ Fourier coefficient via 
\begin{align*}\label{}
\alpha_j(\mu) \approx \eps c_j(\mu)\,.
\end{align*}
One may then hope that an explicit law on the Verblunsky coefficients could shed light on the law of our Fourier coefficients of the GMC. The difficulty here is that the above relationship between the two is only valid in a perturbative way. 

\item The construction of the Virasoro algebra of Liouville conformal field theory (see \cite{BGKRV}), is based on the operators $\mathbf{L}_n$ which generate the Virasoro algebra. Their relationship with the operators $\mathbf{L}_n^0$ in the non-interacting case (the free field case) has the following simple expression:
\begin{align}\label{eqdefLn}
\mathbf{L}_n = \mathbf{L}_n^0 + e^{\gamma c} \int_0^{2\pi}  e^{i n \theta} e^{\gamma \varphi(\theta) d\theta},  \gamma \in (0,2)\,.
\end{align}
As such, our main Theorem says that in the regime $\gamma< \sqrt{2}$, the operators $\mathbf{L}_n$ are asymptotically close (as $n\to \infty$) to the free-field operators $\mathbf{L}_n^0$. This may prove being useful for the analysis of $\mathbf{L}_n$. In the regime $\gamma \in [\sqrt{2}, 2)$ the potential term $e^{\gamma c} \int_0^{2\pi}  e^{i n \theta} e^{\gamma \varphi(\theta) d\theta}$ in \eqref{eqdefLn} is no longer defined via GMC theory but via the theory of Dirichlet forms.

\ei

\vspace{0.1 cm}

Motivated by analytic number theory and random matrix theory, there has been recently some active research around two random sequences which are expected to exhibit the same kind of behaviour as the Fourier coefficients $c_n$ studied in this paper. In the aforementioned models, the n-th term in the sequence implies a finite number of independent variables and hence is closely related to a cut-off version of $c_n$. More specifically, let $\epsilon_n=\frac{1}{n}$ and consider the n-th Fourier coefficient $\tilde{c}_n$ of the measure $e^{\gamma \varphi_{\epsilon_n}(\theta) } d \theta$ which (up to scaling) approximates  $M_\gamma$
\begin{equation*}
\tilde{c}_n= \int_0^{2\pi}  e^{i n \theta} e^{\gamma \varphi_{\epsilon_n}(\theta) } d \theta
\end{equation*}
By the convergence result \eqref{defconvergencechaos}, it is natural to expect that $\tilde{c}_n \approx n^{\frac{\gamma^2}{2}} c_n$ for all $\gamma$. Now, the sequence $\tilde{c}_n$ is expected to display the same behaviour as the following two models:

\begin{itemize}
\item
{\bf Holomorphic multiplicative chaos}: in this setting, one studies the following Fourier coefficients of a complex valued GMC 
\begin{equation*}
a_n= \int_0^{2 \pi} e^{i n \theta} e^{\sqrt{\alpha} X(e^{i\theta})} d \theta
\end{equation*}
where  $\alpha>0$ and $X$ is a holomorphic log-correlated field, i.e. the real part and imaginary part are log-correlated fields which are correlated (in contrast with the case studied in \cite{LRV} where the real part and imaginary part are log-correlated fields which are independent). The field $X$ is defined via the following sum
\begin{equation*}
X(z)= \sum_{k=1}^\infty \frac{X_k}{\sqrt{k}} z^k
\end{equation*}
where $X_k$ are i.i.d. standard complex Gaussian variables. The $a_n$ exhibit the same behaviour as $\tilde{c}_n$ for the value $\alpha=2 \gamma^2$.  In this setting, by exploiting the fact that $X$ is holomorphic, one can see that the coefficients $a_n$ are polynomials in the variables $X_1, \cdots, X_n$ which makes the analysis simpler than the $c_n$ studied in this paper; in particular the $a_n$ have moments of all order. Now, Theorem 1.4 in  \cite{NPS}  is the exact analogue of the capacity estimate  \eqref{capacityestimate} of \cite{LRV} thanks to \eqref{capacityequivalent} (see also \cite{Gerspach} where the author shows the following refined result: for all $\epsilon>0$ the bound $|a_n|= O((\ln n)^{\frac{1}{4}-\epsilon})$ can not hold almost surely for $\gamma= \frac{1}{\sqrt{2}}$). In \cite{NPS}, the authors also prove in Theorem 1.8 the exact analogue of Theorem \ref{thlaw} for the values $\gamma<\frac{1}{2}$; we are able to derive the convergence in law in Theorem \ref{thlaw}  to the optimal bound $\gamma<\frac{1}{\sqrt{2}}$ hence answering Question 1.13 in \cite{NPS} within the framework of this paper, i.e. the study of $c_n$ (for $\gamma \geq \frac{1}{\sqrt{2}}$ a different behaviour is expected).

\item
{\bf Random multiplicative functions}: in this setting one studies partial sums of the form
\begin{equation*}
b_n=\frac{1}{\sqrt{e^{n}}}\sum_{k=1}^{e^n} f(k)
\end{equation*}
where $f$ is a random multiplicative function such that $(f(p))_{\mathcal{P}}$ are i.i.d. where $\mathcal{P}$ denotes the prime numbers. One can consider different distributions for $f(p)$; when $f(p)$ is uniformly distributed on the circle then $b_n$ exhibits the same behaviour as $a_n$ (hence $\tilde{c}_n$) for the value $\gamma= \frac{1}{\sqrt{2}}$ as argued in  \cite{SoundZam} where they prove that $\E[| a_n|]$ is of the order $\frac{1}{(\ln n)^{\frac{1}{4}}}$; the fact that $\E[| b_n|]$ is of the order $\frac{1}{(\ln n)^{\frac{1}{4}}}$ was proved in \cite{Harper1}. In the pioneering work \cite{Harper}, the author shows that the bound  $|b_n|= O((\ln n)^{\frac{1}{4}-\epsilon})$ can not hold almost surely. Finally, the fact that for all $\epsilon>0$ the bound $|b_n|= O((\ln n)^{\frac{1}{4}+\epsilon})$ holds almost surely was proved in \cite{CaicRachid}. The analogue of the results of \cite{Harper,CaicRachid} in our setting is an open problem.

\end{itemize}

\subsection*{Acknowledgments.}
We wish to thank Aleksandr Logunov for an enlightening discussion regarding the relevance of convolutions in Fourier analysis. We thank Reda Chhaibi for a useful discussions about Verblunsky coefficients. 
The first author wishes to thank the university of Geneva far an invited semester in spring 2023 during which this work has been initiated. The second author would like to thank Adam Harper for very interesting discussions on the relation between this paper and related results in other models. The second author would also like to thank Benjamin Bonnefont and Baptiste Cercl\'e for pointing out that Theorem \ref{thlaw} could be proved up to $\gamma^2<\frac{1}{2}$ with the methods of this paper (a previous version of the paper stated the result for $\gamma^2<\frac{1}{3}$). The research of C.G. is supported by the Institut Universitaire de France (IUF), the ERC grant VORTEX 101043450 and the French ANR grant ANR-21-CE40-0003.
The research of V.V. is supported by  the SNSF grant ``2d constructive field theory with exponential interactions''.


\section{Proof of the Rajchman property}

The proof is based on the following simple observation: for $d \geq 1$ an integer, the Fourier coefficient of the $d$-fold convolution $\ast^d M_\gamma$ is $c_n^d$. Since the convergence to $0$ of $c_n$ is equivalent to the convergence to $0$ of $c_n^d$ it is sufficient to study the latter quantity for some suitable choice of $d$. We will in fact show the following: let $d \geq 1$ be some integer, then for $\gamma<\sqrt{\frac{2(d-1)}{d}}$, the $d$-fold convolution $\ast^d M_\gamma$ is a.s. an $L^1$ function. Therefore, by the Riemann-Lebesgue lemma, $c_n^d$ converges to $0$.

Consider the function $f_\epsilon=  \ast^d   \epsilon^{\frac{\gamma^2}{2}}  e^{\gamma \varphi_\epsilon} $ We want to show that there exists $\alpha>0$ such that
\begin{equation}\label{L1conv}
\underset{\epsilon, \epsilon' \to 0}{\lim} \E[ \left ( \int_0^{2 \pi}   |f_\epsilon(u)-f_{\epsilon'}(u)| du \right )^\alpha] =0
\end{equation}
Indeed this will ensure that along a subsequence $f_\epsilon$ converges in $L^1$ which then ensures that the measure $\ast^d M_\gamma$ is an $L^1$ function.

\medskip

The proof is divided into two parts. In order to show \eqref{L1conv}, we will rely on the definition of $f_\epsilon$ as a $d$-fold integral
\begin{equation}\label{deffepsilon}
f_\epsilon(u)=  \int_{[0,2 \pi)^{d-1}} \epsilon^{d \gamma^2/2} e^{\gamma (  \varphi_\epsilon(\theta_1)+ \cdots + \varphi_\epsilon(\theta_{d-1}) +  \varphi_\epsilon(u-\theta_1-\cdots \theta_{d-1})) } d\theta_1 \cdots d \theta_{d-1}
\end{equation}
Then we split the above integral into two pieces, one where the points are close $|\theta_i-\theta_j|\leq \delta$ for some $i,j$ (or $|\theta_i-u+\sum_{j=1}^{d-1} \theta_j|\leq \delta$) and one where the points are all far apart. The first integral can be shown to be small by leveraging capacity estimates obtained on \cite{LRV}. The second integral involves far away points and therefore on this event, the GFFs $\varphi_\epsilon(\theta_{i})$ can be seen as roughly independent of each other, hence rendering the elegant approach from \cite{NathanaelGMC} adaptable.

For $\delta>0$ and any given $u\in[0,2\pi)$, we introduce the following subspace of $[0,2\pi)^{d-1}$ which corresponds to the two-body singularities in the above integral, namely
\begin{align*}\label{}
\mathbf{K}_{\delta,d}(u):=  \bigcup_{-d\leq m \leq d} \left( \bigcup_{1\leq i<j \leq d-1} \{ |\theta_i - \theta_j +2\pi m|  \leq \delta \} \cup  \bigcup_{i=1}^{d-1}\{ |2 \theta_i + \sum_{j\neq i} \theta_j - u +2\pi m| \leq \delta \}  \right)
\end{align*}

Notice that $\mathbf{K}_{\delta, d}(u)$ is made of $d(d-1)/2$ ``slices'' corresponding to $m=0$ plus their $2\pi$ translates. These slices intersect each other and are equal, up to affine change of variables to $\{ |\theta_1-\theta_2|<\delta\}$. We shall often rely on union bounds to reduce the analysis only to one such slice, modulo a multiplicative combinatorial error of $O(d^2)$. We set
\begin{equation*}
F_\eps(u,\theta_1,\ldots,\theta_{d-1}) = \epsilon^{d \gamma^2/2} e^{\gamma (  \varphi_\epsilon(\theta_1)+ \cdots + \varphi_\epsilon(\theta_{d-1}) +  \varphi_\epsilon(u-\theta_1-\cdots \theta_{d-1})) }
\end{equation*}

For any $\alpha\in(0,1)$ (whose value will be fixed later), we have 
\begin{align}\label{e.TooLong}
& \E[ \left ( \int_0^{2 \pi}   |f_\epsilon(u)-f_{\epsilon'}(u)| du \right )^\alpha] \nn \\ & = 
\E[ \left ( \int_0^{2 \pi}   |\int \epsilon^{d \gamma^2/2} e^{\gamma (  \varphi_\epsilon(\theta_1)+ \cdots + \varphi_\epsilon(\theta_{d-1}) +  \varphi_\epsilon(u-\theta_1-\cdots \theta_{d-1})) } - (\eps')^{d \gamma^2/2} \cdots d\theta_1 \cdots d \theta_{d-1}| du \right )^\alpha] \nn \\
& = \Eb{ \left(\int_0^{2\pi}  |\int F_\eps(u,\theta_1,\ldots,\theta_{d-1}) - F_{\eps'}(u,\theta_1,\ldots,\theta_{d-1})  d\theta_1 \cdots d\theta_{d-1} | du \right)^\alpha } \nn \\
& \leq \Eb{ \left(\int_0^{2\pi}  \int_{\mathbf{K}_{\delta,d}(u)} F_\eps(u,\theta_1,\ldots,\theta_{d-1})  d\theta_1 \cdots d\theta_{d-1} du \right)^\alpha }  \nn \\
& + \Eb{ \left(\int_0^{2\pi}  \int_{\mathbf{K}_{\delta,d}(u)} F_{\eps'}(u,\theta_1,\ldots,\theta_{d-1})  d\theta_1 \cdots d\theta_{d-1} du \right)^\alpha } \nn  \\
& + \Eb{ \left(\int_0^{2\pi}  |\int_{\mathbf{K}_{\delta,d}(u)^c} F_\eps(u,\theta_1,\ldots,\theta_{d-1}) - F_{\eps'}(u,\theta_1,\ldots,\theta_{d-1})  d\theta_1 \cdots d\theta_{d-1} | du \right)^\alpha }\,, 
\end{align}
where we used the sub-additivity of $x^\alpha$.

Before we turn to the more problematic third term in~\eqref{e.TooLong}, let us explain why by taking $\alpha$ sufficiently small, one can control the first two terms by a capacity estimate.

\medskip

\ni
\underline{\em The capacity estimate.}
$ $

We now claim that by suitable change of variables, the first term of~\eqref{e.TooLong} can be upper bounded by 
\begin{align*}\label{}
O(d^2)   \E[  \left (  \int_0^{2 \pi}  \int_{|\theta_1-\theta_2| \leq \delta} \epsilon^{d \gamma^2/2} e^{\gamma (  \varphi_\epsilon(\theta_1)+ \cdots + \varphi_\epsilon(u-\theta_1-\cdots \theta_{d-1})) } d\theta_1 \cdots d \theta_{d-1} du \right )^{\alpha}]\,.
\end{align*}
And similarly with $\eps'$-cutoff for the second turn. We now have

 \begin{align*}
&  \E[  \left (  \int_0^{2 \pi}  \int_{|\theta_1-\theta_2| \leq \delta} \epsilon^{d \gamma^2/2} e^{\gamma (  \varphi_\epsilon(\theta_1)+ \cdots + \varphi_\epsilon(u-\theta_1-\cdots \theta_{d-1})) } d\theta_1 \cdots d \theta_{d-1} du \right )^{\alpha}] \\
& \leq  \E \left [ \left ( \int_{|\theta-\theta'| \leq \delta}    M_{\gamma,\epsilon}(d \theta)  M_{\gamma,\epsilon}(d \theta') \right )^{\alpha} M_{\gamma,\epsilon} [0,2 \pi]^{\alpha(d-2)} \right ] \\
& \leq  \E \left [ \left ( \int_{|\theta-\theta'| \leq \delta}    M_{\gamma,\epsilon}(d \theta)  M_{\gamma,\epsilon}(d \theta') \right )^{2 \alpha} \right ]^{1/2} \E[M_{\gamma,\epsilon} [0,2 \pi]^{2 \alpha(d-2)} ]^{1/2}. \\
& \leq \delta^{s \alpha} \E \left [ \left ( \int_0^{2 \pi}   \int_0^{2 \pi}   \frac{M_{\gamma,\epsilon}(d \theta)  M_{\gamma,\epsilon}(d \theta')}{|e^{i\theta} - e^{i\theta'}|^s} \right )^{2 \alpha} \right ]^{1/2} \E[M_{\gamma,\epsilon} [0,2 \pi]^{2 \alpha(d-2)} ]^{1/2} \\
& \leq C \delta^{s \alpha} 
\end{align*}
 where the last inequality is valid for $\alpha$ and $s$ sufficiently small by using \cite{LRV}.

 \medskip
\ni
\underline{\em The $L^1$ convergence}
 $ $

Hence we may now restrict to the third term in~\eqref{e.TooLong}, i.e. to the convolution of points which are all at distance at least $\delta$ from each other. ($\delta$ should be thought as a small but a fixed parameter and $\epsilon,\epsilon' \to 0$). Our first step will be to get rid of the fractional moment $\alpha>0$. Indeed this is needed only to handle the singularities which arise within $\mathbf{K}_{\delta,d}(u)$. We may thus write by using Jensen, 
 
\begin{align*}\label{}
&\Eb{ \left(\int_0^{2\pi}  |\int_{\mathbf{K}_{\delta,d}(u)^c} F_\eps(u,\theta_1,\ldots,\theta_{d-1}) - F_{\eps'}(u,\theta_1,\ldots,\theta_{d-1})  d\theta_1 \cdots d\theta_{d-1} | du \right)^\alpha } \\
& \leq 
\left(\int_0^{2\pi} \Eb{   |\int_{\mathbf{K}_{\delta,d}(u)^c} F_\eps(u,\theta_1,\ldots,\theta_{d-1}) - F_{\eps'}(u,\theta_1,\ldots,\theta_{d-1})  d\theta_1 \cdots d\theta_{d-1} |} \,\,  du  \right)^\alpha\,, 
\end{align*}

Let us state the following Lemma.
\begin{lemma}\label{lemmaL1}
For any $d\geq 2$, any fixed $u\in [0,2\pi)$ and any $\gamma<\gamma_d(:=  \sqrt{2 \frac {d-1} d})$,  along a sufficiently fast subsequence $(\eps_k)_k$ (which does not depend on the choice of $u$ but may depend on the power of convolution $d$), 
\begin{align*}\label{}
\int_{\mathbf{K}_{\delta,d}(u)^c} F_{\eps_k}(u,\theta_1,\ldots,\theta_{d-1})
d\theta_1 \cdots d\theta_{d-1} 
\end{align*}
converges in $L^1$ as $\eps_k \to 0$. 
\end{lemma}

\ni

 
For simplicity let us only stick to the case $u=0$ (the analysis for $u\neq 0$ is very similar as we restrict the integral over $\theta_1,\cdots,\theta_{d-1}$ to $\mathbf{K}_{\delta,d}(u)$ thus avoiding any collapse between the exponentials of the GFF which are involved).  From now on, we will write $\mathbf{K}_{\delta,d}$ for $\mathbf{K}_{\delta,d}(u= 0)$. We may in fact pave $\mathbf{K}_{\delta,d}$ by a finite number of disjoint events which have the following property: there exists intervals $A_i$ for $i= 1 \cdots d$ that are at distance $\delta/2>0$ (we evaluate the distances modulo $2\pi$ here, i.e. using $|\theta-\theta'|_{[2\pi]}:=\inf_{m\in \Z} \{|\theta -\theta'+2\pi m| \} $).
 and such that $\theta_i \in A_i$ for $i= 1 \cdots d-1$ and $-\sum_{j=1}^{d-1} \theta_j \,\, [2\pi]\in A_d $. This is clearly possible and will be used in the sequel. For simplicity we will set $\theta_d=-\sum_{j=1}^{d-1} \theta_j$.

We shall proceed as in \cite{NathanaelGMC}, namely we want to define for any $\alpha>\gamma$ a notion of \textbf{$\alpha$-good points} in such a way that 
\begin{align*}\label{}
\int_{(\theta_1, \cdots, \theta_d) \in A_1 \times \cdots A_d} F_{\eps}(0,\theta_1,\ldots,\theta_{d-1}) 1_{\text{the point }(\theta_1,\ldots,\theta_{d-1})\text{ is $\alpha$-good}}
d\theta_1 \cdots d\theta_{d-1} 
\end{align*}
converges in fact in $L^2$ as $\eps\to 0$, yet containing most of the $L^1$ mass. 
%

The random variable we are interested in is 
\begin{align*}\label{}
\int_{\mathbf{K}_{\delta,d}^c} \eps^{d \gamma^2/2} e^{\gamma \left(\varphi_\eps(\theta_1)+\ldots+\varphi_\eps(\theta_{d-1})+ \varphi_\eps(-(\theta_1 + \ldots + \theta_{d-1})) \right)} d \theta_1 \ldots d \theta_{d-1}
\end{align*}

First of all, Girsanov is telling us what is the suitable notion of typical points. If a point $\Theta:=(\theta_1,\ldots,\theta_{d-1})$ is sampled in $A_1 \times \cdots A_{d-1}$ (and in such a way that $\theta_d\in A_d$) according to the above GMC measure, then conditioned on $\Theta$, by Girsanov, the conditional law of the GFF is given (more or less if we disregard quantities controlled by $\delta$) by a  GFF on $[0,2\pi)$ shifted by the deterministic function
\begin{align*}\label{}
\theta \mapsto \gamma \log \frac 1 {|e^{i\theta_1}  - e^{i \theta}|} + \ldots +
\gamma \log \frac 1 {|e^{i\theta_{d-1}} - e^{i\theta}|} + 
\gamma \log \frac 1 {|e^{i \theta_d} -e^{i \theta}|}  
\end{align*}

This leads us to the following definition of $\alpha$-good $d-1$-tuple of points. 
\begin{definition}\label{}
For any $\alpha>\gamma$ and any $\eps>0$, a $(d-1)$-tuple $\Theta=(\theta_1,\ldots,\theta_{d-1})$ in $\mathbf{K}_{\delta,d}^c$ is called $\alpha$-good at scale $\eps$ if the following event is satisfied
\begin{align*}\label{}
G_\eps^\alpha(\Theta):=  & \{\text{for any $r_1,\ldots, r_{d}$ in $[\eps, \delta]$}, \\ 
& \,\,\,\,\,  \varphi_{r_1}(\theta_1) \leq \alpha \log \frac 1 {r_1},  \ldots, \varphi_{r_{d-1}}(\theta_{d-1}) \leq \alpha \log \frac 1 {r_{d-1}} \text{ and }  
\varphi_{r_d}(\theta_d) \leq \alpha \log \frac 1 {r_d}
\}
\end{align*}
(N.B. Notice that this definition involves  all the scales $r_i$ which are between $\eps$ and the separation distance $\delta$ which controls the punctures). 
\end{definition}

Let us now prove that for any $\gamma<\gamma_d(:=  \sqrt{2 \frac {d-1} d})$, any $\delta>0$ and any $\alpha>\gamma$ sufficiently close to $\gamma$, we have 
\begin{align*}\label{}
\limsup_{\eps \to 0} \Eb{ \eps^{d \gamma^2}  \left(   \int_{ A_1 \times \cdots A_{d-1}} 
e^{\gamma \left(\varphi_\eps(\theta_1)+\ldots+\varphi_\eps(\theta_{d-1})+ \varphi_\eps(-(\theta_1 + \ldots + \theta_{d-1})) \right)}
1_{G_\eps^\alpha(\Theta)} 
d\theta_1 \cdots d\theta_{d-1}   \right)^2}  <\infty 
\end{align*}
In \cite{NathanaelGMC}, proving such an $L^2$ estimate is the main technical step to prove $L^1$ convergence of approximations to a standard GMC: we will thus prove such an $L^2$ estimate and refer to \cite{NathanaelGMC} to promote this estimate to the desired $L^1$ convergence in Lemma \ref{lemmaL1}. We denote the $(d-1)$-tuples of angles:
\begin{align*}\label{}
\Theta:=(\theta_1,\ldots, \theta_{d-1}) \text{    and    } 
\Theta':=(\theta_1',\ldots, \theta_{d-1}') 
\end{align*}
%
On the event $G_\eps^\alpha(\Theta) \cap G_\eps^\alpha(\Theta')$, the following $d$ events are necessarily satisfied:
\begin{align}\label{e.devents}
\begin{cases}
& \mathcal{E}_{r_1}:= \{ \varphi_{r_1}(\theta_1)  \leq  \alpha \log \frac 1 {r_1}\}   \\
&  \ldots\\
& \mathcal{E}_{r_{d-1}}:= \{  \varphi_{r_{d-1}}(\theta_{d-1})  \leq  \alpha \log \frac 1 {r_{d-1}} \}  \\ 
&  \mathcal{E}_{r_d}:= \{   \varphi_{r_d}(-(\theta_1+\ldots+ \theta_{d-1})) \leq  \alpha \log \frac 1 {r_d} \} 
\end{cases}
\end{align}

Under the probability measure
\begin{align*}\label{}
\mathbb{Q}_\eps(d\varphi) := \frac {e^{\gamma 
\left(
\varphi_\eps(\theta_1)+ \ldots + \varphi_\eps(\theta_{d-1})+ \varphi_\eps(-(\theta_1 +  \ldots+ \theta_{d-1}))
+
\varphi_\eps(\theta_1')+\ldots + \varphi_\eps(\theta_{d-1}')+ \varphi_\eps(-(\theta_1' + \ldots+  \theta_{d-1}'))
\right)}
}
{
\Eb{e^{\gamma \left(
\varphi_\eps(\theta_1)+ \ldots + \varphi_\eps(\theta_{d-1})+ \varphi_\eps(-(\theta_1 +  \ldots+ \theta_{d-1}))
+
\varphi_\eps(\theta_1')+\ldots + \varphi_\eps(\theta_{d-1}')+ \varphi_\eps(-(\theta_1' + \ldots+  \theta_{d-1}'))
\right)}}
}\,,
\end{align*}
by Girsanov one has for any $1\leq i \leq d$ (modulo some armless $\delta$-correcting terms), 
\begin{align*}\label{}
& \varphi_{r_i}(\theta_i)  \\
& \sim \mathcal{N}( \Cov{\varphi_{r_i}(\theta_i) , 
\gamma
[
\varphi_\eps(\theta_1)+ \ldots + \varphi_\eps(\theta_{d-1})+ \varphi_\eps(\theta_d)
+ \\
& \hskip 4 cm 
\varphi_\eps(\theta_1')+ \ldots + \varphi_\eps(\theta_{d-1}')+ \varphi_\eps(\theta'_d)
]
},  \log \tfrac 1{r_1})  \\
& \approx \mathcal{N}(  \gamma \log \frac 1 {r_i}+ \gamma \log \frac{1}{r_i+\eps+ |e^{i\theta_i}-e^{i \theta_i'}|},  \log \frac 1 {r_i})\,.
\end{align*}
A crucial observation here is that these $d$ variables are almost uncorrelated from each other  (modulo some $\log \frac 1 \delta$ correlations which we ignore here). This is because we organised the angles in such a way that for each $i$, both $\theta_i$ and $\theta'_i$ are in the same interval $A_i$ and we assumed $A_i$ and $A_j$ to be at least at distance $\delta/2$ when $i\neq j$.
This allows us to obtain the following bound for $r_i=\epsilon+|e^{i\theta_i}-e^{i \theta_i'}|$:
\begin{align*}\label{}
\mathbb{Q}_\eps[\calE_{r_1} \cap \ldots \cap \calE_{r_{d-1}} \cap \calE_{r_d}] 
 & \leq O_{\delta, d}(1) \prod_{i=1}^d (\eps+ |e^{i\theta_i}-e^{i \theta_i'}|)^{\frac{(2 \gamma-\alpha)^2}{2}}
\end{align*}
This leads to
\begin{align*}
& \Eb{ \eps^{d \gamma^2}  \left(   \int_{ A_1 \times \cdots A_{d-1}} 
e^{\gamma \left(\varphi_\eps(\theta_1)+\ldots+\varphi_\eps(\theta_{d-1})+ \varphi_\eps(-(\theta_1 + \ldots + \theta_{d-1})) \right)}
1_{G_\eps^\alpha(\Theta)} 
d\theta_1 \cdots d\theta_{d-1}   \right)^2} \\
& =  \int_{ A_1 \times \cdots A_{d-1} \times  A_1 \times \cdots A_{d-1} } 
 \eps^{d \gamma^2}  \Eb{ e^{\gamma \left(
\varphi_\eps(\theta_1)+\varphi_\eps(\theta_{2})+ 
\ldots 
 \right)}
1_{G_\eps^\alpha(\Theta)}
1_{G_\eps^\alpha(\Theta')} }
d\Theta d \Theta' \\
& \leq 
\int_{ A_1 \times \cdots A_{d-1} \times  A_1 \times \cdots A_{d-1} }  \eps^{d \gamma^2}  
\Eb{ e^{\gamma \left(
\varphi_\eps(\theta_1)+\varphi_\eps(\theta_{2})+ 
\ldots 
 \right)}}
\mathbb{Q}_\eps[\calE_{r_1} \cap \ldots \cap \calE_{r_{d-1}} \cap \calE_{r_d}] 
 d\Theta d \Theta' 
 \\
& \leq
O_{\delta, d}(1) 
\int_{ (A_1 \times \cdots A_{d-1})^2 } 
 \prod_{i=1}^d  \frac 1 {(\eps+ |e^{i\theta_i}-e^{i \theta_i'}|)^{\gamma^2}} 
 \prod_{i=1}^d (\eps+ |e^{i\theta_i}-e^{i \theta_i'}|)^{\frac{(2 \gamma-\alpha)^2}{2}} 
 d\Theta d \Theta' 
 \\
& \leq
O_{\delta, d}(1) 
\int_{ (A_1 \times \cdots A_{d-1})^2 } 
 \prod_{i=1}^d  \frac 1 {(\eps+ |e^{i\theta_i}-e^{i \theta_i'}|)^{\gamma^2  - \frac{(2 \gamma-\alpha)^2}{2}}} 
 d\Theta d \Theta' \,.
\end{align*}
%

Now let us choose $\alpha>\gamma$ close enough to $\gamma$  so that the exponent above satisfies 
\begin{align*}\label{}
\gamma^2  - \frac{(2 \gamma-\alpha)^2}{2} < \frac {d-1} d
\end{align*}
(Recall we assumed $\gamma^2 < 2 \frac {d-1} d$). To conclude the proof of Theorem \ref{Rajchmanback}, note that by applying in the above integral the following change of variables: $x_i:=\theta_i-\theta'_i$, $y_i=\theta_i+\theta_i'$, and using the integral estimate below from Lemma \ref{l.integral}, we obtain that the above $L^2$ estimate indeed remains finite as $\eps\searrow 0$, which thus  concludes the proof.
\qed

\medskip
 
\begin{lemma}\label{l.integral}
Let $d \geq 1$ be an integer. For all $A>0$ and $u < \frac{d-1}{d}$, the following holds
\begin{equation}\label{intconv}
\int_{ [-A,A]^{d-1}} \left (  \prod_{i=1}^{d-1} \frac{1}{|x_i|^{u}}  \right ) \frac{1}{|\sum_{i=1}^{d-1} x_i|^{u}}  dx_1 \cdots d x_{d-1}< \infty
\end{equation}
\end{lemma}

\proof

First  suppose that $|x_1| <  \cdots < |x_{d-1}| <  |\sum_{i=1}^{d-1} x_i|$. Then we have
\begin{align*}
& \int_{ |x_1| <  \cdots < |x_{d-1}| <  |\sum_{i=1}^{d-1} x_i|} \left (  \prod_{i=1}^{d-1} \frac{1}{|x_i|^{u}}  \right ) \frac{1}{|\sum_{i=1}^{d-1} x_i|^{u}}  dx_1 \cdots d x_{d-1} \\
& \leq \int_{ |x_1| <  \cdots < |x_{d-1}| }  \frac{1}{|x_{1}|^{u}}  \cdots  \frac{1}{|x_{d-2}|^{u}}  \frac{1}{|x_{d-1}|^{2u}}  dx_1 \cdots d x_{d-1} \\ 
& \leq C \int_{ |x_1| <  \cdots < |x_{d-2}| }  \frac{1}{|x_{1}|^{u}}  \cdots  \frac{1}{|x_{d-2}|^{3u-1}}   dx_1 \cdots d x_{d-2} \\ 
& \leq C \int_{ |x_1| <A }  \frac{1}{|x_{1}|^{d u-(d-2) }}    dx_1  \\ 
& \leq C \\
\end{align*}

Now suppose that the maximum of $|x_1|,  \cdots,  |x_{d-1}|,  |\sum_{i=1}^{d-1} x_i|$ is achieved say at $x_1$. Then we perform the change of variable $x_1=  \sum_{i=1}^{d-1} v_i$ and $x_2= -v_2, \cdots x_{d-1}= -v_{d-1}$, in which case $\sum_{i=1}^{d-1} x_i= v_1$ and we are back to the previous case.


\qed

\section{Proof of Theorem \ref{RajchmanbackL4} via $L^4$ moments} 

\noindent
{\em Proof of Theorem \ref{RajchmanbackL4}.}

By a simple computation on integrals and using Fubini, we get
\begin{align*}
& \E[  |  c_n   |^4  ]   \\
& =\int_{[-\pi, \pi]^4}   e^{i n (\theta_1-\theta_2 +\theta_3-\theta_4)}\prod_{j<k} \frac{1}{|e^{i \theta_j}-e^{i \theta_k} |^{\gamma^2}} d\theta_1 d\theta_2 d\theta_3 d\theta_4 \\
& \leq C 
\int_{[-\frac \pi 2, \frac \pi 2]^4}   e^{i n (\theta_1-\theta_2 +\theta_3-\theta_4)}\prod_{j<k} \frac{1}{|e^{i \theta_j}-e^{i \theta_k} |^{\gamma^2}} d\theta_1 d\theta_2 d\theta_3 d\theta_4\,.
\end{align*}
By reducing the analysis to the smaller cube $[-\frac \pi 2, \frac \pi 2]^4$, this allows us to use the lower bound: 
\begin{align*}\label{}
|e^{i \theta_j} -e^{i \theta_k}| \geq c |\theta_j - \theta_k|
\end{align*}
for any $j<k$ for some small enough constant $c>0$. This gives us 

\begin{align*}
& \E[  |  c_n   |^4  ] \\
& \leq \frac{O(1)}{n^{4-6\gamma^2}}
\int_{[-n \tfrac \pi 2, n \tfrac \pi 2]}   e^{i  (\theta_1-\theta_2+\theta_3-\theta_4)}\prod_{j<k} \frac{1}{|\theta_j-\theta_k |^{\gamma^2}} d\theta_1 d\theta_2 d\theta_3 d\theta_4 \\
& = \frac{O(1)}{n^{4-6\gamma^2}}\int_{-n \tfrac \pi 2}^{n \tfrac \pi 2} \int_{-n \tfrac \pi 2}^{n \tfrac \pi 2}   \int_{-n \tfrac \pi 2}^{n\tfrac \pi 2}   \frac{e^{i(\theta_3-\theta_4)}}{\prod_{2 \leq i<j}|\theta_i- \theta_j|^{\gamma^2}}   \left ( \int_{-n\frac \pi 2}^{n \tfrac \pi 2}  \frac{e^{i(\theta_1-\theta_2)}}{\prod_{1<j}|\theta_1 - \theta_j|^{\gamma^2}}  d\theta_1 \right )  d\theta_2 d\theta_3 d\theta_4  \\
& \leq   \frac{O(1)}{n^{4-6\gamma^2}}\int_{-n\tfrac \pi 2 }^{n \tfrac \pi 2} \int_{-n \tfrac \pi 2}^{n \tfrac \pi 2}   \int_{-n \tfrac \pi 2}^{n \tfrac \pi 2}   \frac{e^{i(\theta_3-\theta_4)}}{\prod_{2 \leq i<j}|\theta_i - \theta_j|^{\gamma^2}}  F(  \theta_2,\theta_3,\theta_4  )  d\theta_2 d\theta_3 d\theta_4\,,   \\
\end{align*}
where 
\begin{equation*}
F(  \theta_2,\theta_3,\theta_4  )= \int_\R  \frac{e^{iv}}{  |v|^{\gamma^2} |v+ \theta_2-\theta_3|^{\gamma^2}  |v+ \theta_2-\theta_4|^{\gamma^2}  }   dv \,. 
 \end{equation*}
Then integrating with respect to $\theta_3$ we get
\begin{align*}
& \E[  |  a_n   |^4  ]   \\
& \leq  \frac{O(1)}{n^{4-6\gamma^2}}\int_{-n \tfrac \pi 2}^{n\tfrac \pi 2} \int_{-n\tfrac \pi 2}^{n \tfrac \pi 2}   \int_{-n \tfrac \pi 2}^{n \tfrac \pi 2}   \frac{e^{i(\theta_3-\theta_4)}}{\prod_{2 \leq i<j}|\theta_i-\theta_j|^{\gamma^2}}  F(  \theta_2,\theta_3,\theta_4  )  d\theta_2 d\theta_3 d\theta_4  \\
& \leq \frac{O(1)}{n^{4-6\gamma^2}}\int_{-n\tfrac \pi 2}^{n \tfrac \pi 2} \int_{-n \tfrac \pi 2}^{n \tfrac \pi 2}   \frac{1}{|\theta_2-\theta_4|^{\gamma^2}}  \left ( \int_{-\infty}^{\infty}   \frac{e^{i w }}{ |w|^{\gamma^2}  |w+\theta_4-\theta_2|^{\gamma^2}  }  F(  \theta_2,w+\theta_4,\theta_4  ) dw \right )  d\theta_2 d\theta_4  \\
\end{align*}
Now, when $|\theta_2-\theta_4|$ is large then
\begin{equation*}
\int_{-\infty}^{\infty}   \frac{e^{i w }}{ |w|^{\gamma^2}  |w+\theta_4-\theta_2|^{\gamma^2}  }  F(  \theta_2,w+\theta_4,\theta_4  ) dw \sim \frac{\kappa^2}{|\theta_2-\theta_4|^{3 \gamma^2}}
\end{equation*}
hence one has
\begin{align*}
  \E[  |  a_n   |^4  ]   & \leq  \frac{O(1)}{n^{4-6\gamma^2}}\int_{-n\tfrac \pi 2}^{n \tfrac \pi 2} \int_{-n \tfrac \pi 2}^{n \tfrac \pi 2}   \frac{1}{|\theta_2-\theta_4|^{4 \gamma^2}}   d\theta_2 d\theta_4  \\
& \leq \frac{O(1)}{ n^{2-2\gamma^2}  }    \int_{-\tfrac \pi 2}^{\tfrac \pi 2} \int_{- \tfrac \pi 2}^{\tfrac \pi 2}   \frac{1}{|v-v'|^{4 \gamma^2}}   dv dv' \,. \\
\end{align*}
When $\gamma^2 < \frac 1 2 $, this leads us to the upper bound uniformly on $n\geq 1$, 
\begin{align*}
\E[|c_n|^4] & \leq  \frac{C}{n^{2(1-\gamma^2)}}. \\
\end{align*}

Therefore if $\beta<\frac{1-2 \gamma^2}{4}$, one gets that $\E [\sum_{n\geq 1}(|c_n|n^\beta)^4]< \infty$ and in particular almost surely $|c_n|n^\beta$ converges to $0$ when $n$ goes to infinity.

\qed

\begin{remark}
The method of moments to prove Theorem \ref{RajchmanbackL4} can easily be adapted to the case of \cite{FalkJin} where it is likely to produce more optimal bounds than their Corollary 2.6 on the Fourier decay of the 2d GMC measure.
\end{remark}

\vspace{0.1 cm}

\section{Convergence in law (Theorem \ref{thlaw})}

Our second main theorem (Theorem \ref{thlaw}) deals with the convergence in law of the (suitably rescaled) Fourier coefficients
\begin{equation*}
c_n=c_n(M_\gamma) =\int_0^{2\pi}e^{\gamma \varphi(\theta)-\frac{\gamma^2}{2}\E[\varphi(\theta)^2]} e^{i n \theta} d\theta 
\end{equation*}
where we recall that  $\varphi$ is the GFF on the circle $\S^1$ with covariance
\begin{equation*}
\E[\varphi(\theta) \varphi(\theta')]= \ln \frac{1}{|e^{i \theta}- e^{i \theta'}|} 
\end{equation*}

The proof of Theorem \ref{thlaw} will rely on the independence properties of a certain ``white-noise'' decomposition of the Gaussian field $\varphi$ on the circle $\S^1$ as constructed in \cite{junnila2017uniqueness}. This white-noise decomposition, though very nice, is slightly frightening at first sight. To make the main ideas more transparent we will  prove a convergence in law for a log-correlated field on $[0,1]=\R / \Z$ whose white-noise decomposition is easier to handle and goes back to \cite{bacry2003log}. We will then only shortly explain in Section \ref{ss.circle} how to adapt the proof to the case of the circle.

\subsection{Convergence in law for a ``toy-model'' on $[0,1]$.}
$ $

We consider the following Gaussian log-correlated field $X$ on the unit interval $[0,1]$ (where the positions are not viewed modulo 1):
\begin{equation*}
\E[X(x) X(y)]= \ln \frac{1}{|x-y|}
\end{equation*}
(in particular the distance $|x-y|$ is not computed modulo 1). 

For each $n\in \Z$, we let $Z_n$ be the $n^{th}$ Fourier coefficient of the $\gamma$-GMC induced by $X$ and viewed as a $1$-periodic measure. Namely,
\begin{align*}\label{}
Z_n:= \int_0^1 e^{\gamma X(x) - \frac{\gamma^2} 2 \Eb{X(x)^2}} e^{2 i \pi n x} dx\,.
\end{align*}

 We introduce $\sigma$ via the following formula
\begin{equation}\label{e.SIGMA}
\sigma^2 = \sigma^2(X):=  \left ( \int_{\R}  \frac{e^{i 2\pi u }}{|u|^{\gamma^2}}  du   \right )  \int_{[0, 1] } e^{2 \gamma X(x) - 2 \gamma^2 \E[  X(x) ] }   dx 
\end{equation}
Our  main theorem on the  convergence in law in this ``toy-model'' case reads as follows:
\begin{theorem}\label{th.Toy}
For any $\gamma^2 <\frac{1}{3}$, the complex random variable $n^{\frac{1-\gamma^2}{2}} Z_n$ converges in law to $\frac \sigma {\sqrt{2}} N$ where $N$ is a standard complex variable independent from $\sigma$.
\end{theorem}

We will rely extensively in the proof on the following important fact: as observed in \cite{bacry2003log}, there exists a very convenient white-noise decomposition of the field $X$ which  for any $t>1$, induces a ``$t$-regularization'' $X_t$ of $X$ with explicit covariance structure 
\begin{equation}\label{e.BB}
\begin{cases}
& \E[X_t(x)X_t(y)]= \ln \frac{1}{|x-y|} \quad \text{if} \quad |x-y|>\frac{1}{t} \\
& \E[X_t(x)X_t(y)]= \ln t +1 -t |x-y| \quad \text{if} \quad |x-y| \leq \frac{1}{t}\,.
\end{cases}
\end{equation}
Furthermore, $X_t$ is a random field with independent increments (as $t$ grows).

\medskip
\ni
{\em Proof of Theorem \ref{th.Toy}.} 

Let $A>0$ be fixed and let us write 
\begin{equation*}
Z_n=    Z_n- \E[ Z_n |  \mathcal{F}_{\frac{n}{A}} ]+\E[ Z_n |  \mathcal{F}_{\frac{n}{A}} ]\,.
\end{equation*}
We have the following 
\begin{equation*}
\E[ Z_n |  \mathcal{F}_{\frac{n}{A}} ] =\int_0^1e^{\gamma X_{\frac{n}{A}}(x)-\frac{\gamma^2}{2}\E[X_{\frac{n}{A}}(x)^2]} e^{i 2\pi n x} dx
\end{equation*}
(N.B. This martingale property is a very convenient feature of this particular regularisation of the log-correlated field $X$).
Furthemore with $K_{\frac{n} A}(|x-y|):= \Eb{X(x)X(y)} - \Eb{X_{n/A}(x) X_{n/A}(y)}$, i.e.  
\begin{equation*}
\begin{cases}
& K_{\frac{n}{A} } (|x-y|) := \ln \frac{A}{n |x-y|} -1  + \frac{n |x-y|}{A}\quad \text{if} \quad |x-y| \leq \frac{A}{n} \\
&  K_{\frac{n}{A} } (|x-y|) :=0 \quad \text{otherwise}
\end{cases}
\end{equation*}
we get

\begin{align*}
& \E[|Z_n- \E[ Z_n |  \mathcal{F}_{\frac{n}{A}} ] |^2  |  \mathcal{F}_{\frac{n}{A}}  ] \\
& = \int_{[0,1]^2} e^{\gamma X_{\frac{n}{A}}(x)+\gamma X_{\frac{n}{A}}(y) - \gamma^2 (\ln \frac{n}{A} +1 )} e^{i 2\pi n (x-y)}  \left ( e^{\gamma^2 K_{\frac{n}{A} } (|x-y|)} -1 \right ) dx dy  \\
& = \int_{x\in[0,1], |x-y|\leq A/n} e^{ \gamma X_{\frac{n}{A}}(x) +
\gamma X_{\frac{n}{A}}(y) - \gamma^2 (\ln \frac{n}{A} +1 )} 
e^{i 2\pi n (x-y)}  \left ( e^{\gamma^2 K_{\frac{n}{A} } (|x-y|)} -1 \right ) dx dy \\
& =   (n/A)^{\gamma^2} e^{\gamma^2} \frac 1 n \int_{x\in [0,1], |u|\leq A } 
e^{ \gamma X_{\frac{n}{A}}(x) +
\gamma X_{\frac{n}{A}}(x+\frac 1 n u) - 2 \gamma^2 (\ln \frac{n}{A} +1 )}  \\
& \hskip 6 cm \times 
e^{2i \pi u }
 \left ( \frac{A^{\gamma^2}e^{-\gamma^2 + \gamma^2 \frac{|u|}{A}}}{|u|^{\gamma^2}} -1 \right ) dx du  \\
& = n^{\gamma^2-1}\int_{x\in [0,1], |u|\leq A } 
e^{ \gamma X_{\frac{n}{A}}(x) +
\gamma X_{\frac{n}{A}}(x+\frac 1 n u) - 2\gamma^2 (\ln \frac{n}{A} +1 )} 
e^{2i \pi u }
 \left ( \frac{e^{ \gamma^2 \frac{|u|}{A}}}{|u|^{\gamma^2}} - \frac{e^{\gamma^2}}{A^{\gamma^2}} \right ) dx du  \\
& = n^{\gamma^2-1}\int_{x\in [0,1], |u|\leq A } 
 e^{ \gamma (X_{\frac{n}{A}}(x)+  X_{\frac{n}{A}}(x+ \frac {1}{n} u)) - \frac{\gamma^2} {2}  \E[(X_{\frac{n}{A}}(x)+  X_{\frac{n}{A}}(x+ \frac {1}{n} u))^2] }  \\
& \hskip 6 cm 
e^{2i \pi u }
 \left ( \frac{1}{|u|^{\gamma^2}} -\frac{ e^{\gamma^2} e^{-\gamma^2 \frac {|u|}{A}} }{A^{\gamma^2}}  \right )\,,
 dx du   
\end{align*}
where in the last line, we absorbed the term $e^{\gamma^2 \frac {|u|}{A}}$ in order to renormalize by the variance $ \frac{\gamma^2} 2 \E[(X_{\frac{n}{A}}(x)+  X_{\frac{n}{A}}(x+ \frac {1}{n} u))^2]$ instead of $2 \gamma^2 (\ln \frac{n}{A} +1 )$.


From now on, we set 
\begin{equation*}
\begin{cases}
& \sigma_{A,n}^2 := n^{1-\gamma^2} \E[ \left|Z_n- \E[ Z_n |  \mathcal{F}_{\frac{n}{A}} ] \right|^2  |  \mathcal{F}_{\frac{n}{A}}  ] \\
& \sigma_{A,n, \R}^2 := n^{1-\gamma^2} \E[\Re \left(Z_n- \E[ Z_n |  \mathcal{F}_{\frac{n}{A}} ] \right)^2  |  \mathcal{F}_{\frac{n}{A}}  ] \\
& \sigma_{A,n,  i \R}^2 := n^{1-\gamma^2} \E[\Im \left(Z_n- \E[ Z_n |  \mathcal{F}_{\frac{n}{A}} ] \right)^2  |  \mathcal{F}_{\frac{n}{A}}  ]  
\end{cases}
\end{equation*}

Now, we have the following lemma
\begin{lemma}\label{l.simpler}
The coupled random variables $(\sigma_{A,n, \R}^2, \sigma_{A,n, i\R}^2)$  converge in law as $n$ goes to infinity to zero, to the random  variable
\begin{align*}\label{}
\frac 1 2 (\sigma_A^2, \sigma_A^2)\,,
\end{align*}
where
\begin{equation*}
\sigma_A^2=  \left ( \int_{|u| \leq A}   e^{i  2\pi u }  \left ( \frac{1}{|u|^{\gamma^2}} -\frac{ e^{\gamma^2} e^{-\gamma^2 \frac {|u|}{A}} }{A^{\gamma^2}} \right ) du \right )  \int_{[0,1] } e^{2 \gamma X(x) - 2 \gamma^2 \E[  X(x) ] } dx 
\end{equation*}
\end{lemma}

\begin{remark}\label{}
We expect that a convergence in probability should hold here.  This may follow by applying the techniques from \cite{shamov2016gaussian} (while \cite{robert2010gaussian} focuses on the uniqueness for the convergence in law). 
\end{remark}

\proof

One writes 
\begin{align*}
 & \sigma_{A,n}^2 \\
 & = \int_{[0,1], |u| \leq A } e^{ \gamma (X_{\frac{n}{A}}(x)+  X_{\frac{n}{A}}(x+ \frac {1}{n} u)) -  \frac{\gamma^2} 2  \E[(X_{\frac{n}{A}}(x)+  X_{\frac{n}{A}}(x+ \frac {1}{n} u))^2] }  e^{i 2 \pi u}   \left ( \frac{ 1}{|u|^{\gamma^2}} -\frac{e^{\gamma^2} e^{-\gamma^2 \frac {|u|}{A}}}{A^{\gamma^2}} \right ) dx du 
\end{align*}
Now we have the following property: the quantity  
\begin{equation*}
\limsup_{n\to \infty} \sup_{|x-y|>\frac{C}{n}, |u|, |v| \leq A}| \E[ (X_{\frac{n}{A}}(x)+  X_{\frac{n}{A}}(x+ \frac {1}{n} u) (X_{\frac{n}{A}}(y)+  X_{\frac{n}{A}}(y+ \frac {1}{n} v) ]  - 4\E[ (X_{\frac{n}{A}}(x) X_{\frac{n}{A}}(y) ]   |
\end{equation*}

goes to $0$ as $C$ goes to infinity and the quantity 
\begin{equation*}
\limsup_{n\to \infty} \sup_{x,y \in [0,1]^2, |u|, |v| \leq A}| \E[ (X_{\frac{n}{A}}(x)+  X_{\frac{n}{A}}(x+ \frac {1}{n} u) (X_{\frac{n}{A}}(y)+  X_{\frac{n}{A}}(y+ \frac {1}{n} v) ]  - 4\E[ (X_{\frac{n}{A}}(x) X_{\frac{n}{A}}(y) ]   |
\end{equation*}
is bounded. The desired result  can then be obtained by an adaptation of the work \cite{robert2010gaussian} by Robert and the second author. See in particular Lemma 4.3 as well as the arguments below which prove that GMC measures do not depend on the chosen mollifier kernel $\theta$.

\medskip
%

\medskip

Now let us discuss briefly how to handle $ \sigma_{A,n,\R}^2$ and $ \sigma_{A,n, i\R}^2$. The idea is to show that the fluctuations are asymptotically equally split on the real and imaginary parts. To see this, let us write 
\begin{align*}\label{}
\sigma_{A,n,\R}^2 & = 
n^{1-\gamma^2} 
 \int_{[0,1]^2} e^{\gamma X_{\frac{n}{A}}(x)+\gamma X_{\frac{n}{A}}(y) - \gamma^2 (\ln \frac{n}{A} +1 )} \cos(2\pi n x) \cos(2\pi ny) \\
 & \hskip 6cm \times   \left ( e^{\gamma^2 K_{\frac{n}{A} } (|x-y|)} -1 \right ) dx dy\,. 
\end{align*}
Now using the identity 
\begin{align*}\label{}
\cos(2\pi n x) \cos(2\pi ny) = \frac 1 2 \cos(2\pi n (x-y)) + \frac 1 2 \cos(2\pi n (x+y))\,,
\end{align*}
we readily see that the integral corresponding to the part $\frac 1 2 \cos(2\pi n (x-y))$ will converge in law as in the above case to $\frac 1 2  \sigma_A^2$, while the second part corresponding to $ \frac 1 2 \cos(2\pi n (x+y))$ will converge instead in probability to zero because the high oscillating mode $n$ is not cancelled out by the scaling window $|x-y|\leq \frac A n$.

%
%
\qed

\medskip
Now, if $X=X_{\frac{n}{A}}+X^{>}_{\frac{n}{A}}$ is the decomposition of $X$ in two independent fields (recall that $K_{\frac{n}{A} } $ is the covariance of $X^{>}_{\frac{n}{A}}$ ), by Theorem 2.6 in \cite{ChenShao}, applied to
\begin{equation*}
Y_{i,\R} := \frac{n^{\frac{1-\gamma^2}{2}}}{\sigma_{A,n,\R}} \int_{\frac{(i-1)A}{n}}^{\frac{ iA}{n}} e^{\gamma X_{\frac{n}{A}}(x)-\frac{\gamma^2}{2}\E[X_{\frac{n}{A}}(x)^2]} \cos(2 \pi nx) ( e^{\gamma X^{>}_{\frac{n}{A}}(x)-\frac{\gamma^2}{2}\E[X^{>}_{\frac{n}{A}}(x)^2]}-1 ) dx 
\end{equation*}
under the measure $\P(. |  \mathcal{F}_{\frac{n}{A}} )$, for all $\eta \in (0,1)$ we get the existence of an absolute constant $C>0$ such that
\begin{equation*}
\sup_{z \in \R} | \P( \Re(Z_n- \E[ Z_n |  \mathcal{F}_{\frac{n}{A}} ]) \leq z \,  n^{\frac{\gamma^2 -1}{2}} |  \mathcal{F}_{\frac{n}{A}} )  - \P( \sigma_{A,n,\R} N  \leq z  |  \mathcal{F}_{\frac{n}{A}} )   | \leq C \sum_{i=1}^{\frac{n}{A}} E[ |Y_{i,\R}|^{2+\eta} |  \mathcal{F}_{\frac{n}{A}}] 
\end{equation*}
where $N$ is a Gaussian independent from everything. Let us fix $z\not =0$, say $z<0$ (the other case $z>0$ can be dealt similarly). For all $\epsilon>0$, we have 
\begin{align*}
&  | \P( \Re (Z_n- \E[ Z_n |  \mathcal{F}_{\frac{n}{A}} ] ) \leq z \,  n^{\frac{\gamma^2 -1}{2}} )  - \P( \sigma_{A,n,\R} N  \leq z   )   |    \\
& \leq  \E[   |     \P(  \Re(Z_n- \E[ Z_n |  \mathcal{F}_{\frac{n}{A}} ] ) \leq z \,  n^{\frac{\gamma^2 -1}{2}}  |  \mathcal{F}_{\frac{n}{A}} )  - \P( \sigma_{A,n,\R} N  \leq z  |  \mathcal{F}_{\frac{n}{A}} )     |   ]  \\
& \leq  \E[   |     \P( \Re(  Z_n- \E[ Z_n |  \mathcal{F}_{\frac{n}{A}} ] ) \leq z \,  n^{\frac{\gamma^2 -1}{2}} |  \mathcal{F}_{\frac{n}{A}} )  - \P( \sigma_{A,n,\R} N  \leq z  |  \mathcal{F}_{\frac{n}{A}} )  |1_{\sigma_{A,n,\R} \leq \epsilon}        ]  \\
& + \E[   |     \P( \Re(Z_n- \E[ Z_n |  \mathcal{F}_{\frac{n}{A}} ] )\leq z \,  n^{\frac{\gamma^2 -1}{2}} |  \mathcal{F}_{\frac{n}{A}} )  - \P( \sigma_{A,n,\R} N  \leq z  |  \mathcal{F}_{\frac{n}{A}} ) | 1_{\sigma_{A,n,\R}> \epsilon}   ]  \\
& \leq 2 \frac{\epsilon^2}{z^2}  + \frac{C}{\epsilon^{2+\eta}} \sum_{i=1}^n \E[ \bar{Y}_i^{2+\eta}]        
\end{align*}
where we have used Markov to bound the first term and 
\begin{equation*}
\bar{Y}_i := n^{\frac{1-\gamma^2}{2}} \int_{\frac{(i-1)A}{n}}^{\frac{ iA}{n}} e^{\gamma X_{\frac{n}{A}}(x)-\frac{\gamma^2}{2}\E[X_{\frac{n}{A}}(x)^2]}  ( e^{\gamma X^{>}_{\frac{n}{A}}(x)-\frac{\gamma^2}{2}\E[X^{>}_{\frac{n}{A}}(x)^2]}+1 ) dx 
\end{equation*}
Recall the classical estimate (see \cite{RV}) 
\begin{equation*}
E[ \bar{Y}_i^{2+\eta} ] \leq C \frac{n^{\frac{(2+\eta)(1-\gamma^2)}{2}}}{n^{\zeta(2+\eta)}}
\end{equation*}
where $\zeta(q)=(1+\frac{\gamma^2}{2})q-\frac{\gamma^2}{2}q^2 $. Therefore $n E[ \bar{Y}_i^3 ]$ converges to $0$ for $\gamma^2< \frac{1}{2}$ and $\eta$ sufficiently small since $1+ \frac{(2+\eta)(1-\gamma^2)}{2} - \zeta(2+\eta)= \eta(\gamma^2-\frac{1}{2})+O(\eta^2)$. In conclusion, for all $z \not = 0$ and $\epsilon>0$
\begin{equation*}
 \underset{n \to \infty} {\overline \lim} | \P( \Re( Z_n- \E[ Z_n |  \mathcal{F}_{\frac{n}{A}} ]) \leq z  n^{\frac{\gamma^2-1} 2}  )  - \P( \sigma_{A,n,\R} N  \leq z   ) |  \leq 2 \frac{\epsilon^2}{z^2} 
\end{equation*}
and therefore $n^{\frac{1-\gamma^2} 2}\, \Re(Z_n- \E[ Z_n |  \mathcal{F}_{\frac{n}{A}} ])$ converges in law as $n$ goes to infinity to $ \frac{\sigma_A}{\sqrt{2}} N  $.

Finally, we can conclude for the convergence of $n^{\frac{1-\gamma^2} 2} Re(Z_n)$ by registering the two following facts:
\begin{align*}
 (1) \;\;\;\;  & \sigma_A^2 \text{  converges a.s. as $A \to \infty$ to } \sigma^2 \\
 (2) \;\;\;\; & \underset{n \to \infty}{\overline \lim} \:   n^{1-\gamma^2}\E[  | \E[ Z_n |  \mathcal{F}_{\frac{n}{A}} ]|^2  ] \underset{A \to \infty}{ \rightarrow} 0 
\end{align*}
The first fact is just a consequence of the deterministic convergence 
\begin{align*}\label{}
\int_{|u| \leq A}   e^{i  2\pi u }  \left ( \frac{1}{|u|^{\gamma^2}} -\frac{ e^{\gamma^2} e^{-\gamma^2 \frac {|u|}{A}} }{A^{\gamma^2}} \right ) du   \underset{A\to\infty}\longrightarrow
\kappa(\gamma) = \int_{\R}   e^{i  2\pi u }   \frac{1}{|u|^{\gamma^2}}  du \,.
\end{align*}

The second fact, i.e. the convergence to $0$ is not straightforward and is the object of the following lemma:

\begin{lemma}\label{}
\begin{align*}\label{}
\limsup  n^{1-\gamma^2} \Eb{ | \Eb{Z_n \md  \calF_{n/A}} |^2}  \underset{A\to \infty}\longrightarrow  0
\end{align*}
\end{lemma}

\proof

We have
\begin{align*}\label{}
n^{1-\gamma^2}\Eb{ | \Eb{Z_n \md  \calF_{n/A}} |^2} & = n^{1-\gamma^2} \Eb{\int_{[0,1]^2} e^{i n 2\pi (y-x)}  e^{\gamma (X_{n/A}(x)+ X_{n/A}(y)) - \gamma^2 (\log \frac n A +1) } dxdy} \\
& = n^{1-\gamma^2} \int_{|x-y| \leq \frac{A}{n}} e^{i n 2 \pi (y-x)} e^{\gamma^2 ( \ln \frac{n}{A}+ 1-\frac{n}{A}|x-y| ) } dx dy \\
& + n^{1-\gamma^2} \int_{|x-y| > \frac{A}{n}} e^{i n 2 \pi (y-x)} \frac{1}{|x-y|^{\gamma^2}}dx dy    \\
& = \frac{1}{A^{\gamma^2}} \int_{|u| \leq A} e^{i 2 \pi u} e^{\gamma^2(1- \frac{|u|}{A})} du+  \int_{|u| >  A}  \frac{e^{i u }}{|u|^{\gamma^2}}du
\end{align*}
which converges to $0$ as $A$ goes to infinity.

\medskip
\medskip
We obtain in the exact same fashion the convergence in law of the imaginary part: 
\begin{align*}\label{}
n^{\frac{1-\gamma^2} 2} \Im(Z_n)
\end{align*}

The only property, we still need to justify is the fact that  $(n^{\frac{1-\gamma^2} 2} \Re(Z_n)), n^{\frac{1-\gamma^2} 2} \Im(Z_n))$ is asymptotically (in law) a Gaussian vector under the conditional measure $P( \cdot |  \mathcal{F}_{\frac{n}{A}} )$ with covariance matrix $\frac {\sigma^2} 2 I_2$ (where $I_2$ is the identity matrix on $\R^2$). 
To show that the limiting vector is Gaussian, one can run the same analysis, i.e. under $P( \cdot |  \mathcal{F}_{\frac{n}{A}} )$, on linear combinations $\lambda_{\R} n^{\frac{1-\gamma^2} 2} \Re(Z_n) + \lambda_{i \R} n^{\frac{1-\gamma^2} 2} \Im(Z_n)$. And finally, to identify the correct covariance matrix, we are left with proving that in law,
\begin{align*}\label{}
n^{1-\gamma^2} \Eb{\Re(Z_n) \Im(Z_n) \md  \calF_{\frac n A} } \to 0\,.
\end{align*}
As above, instead of working with $Z_n$, the main step is to work instead with $Z_n - \Eb{Z_n \md \calF_{\frac n A}}$ and to identify that, in law,

\begin{align*}\label{}
n^{1-\gamma^2}
\Eb{\Re(Z_n \, -  \Eb{Z_n \md \calF_{\frac n A}}) 
\Im(Z_n \, -  \Eb{Z_n \md \calF_{\frac n A}} \md \calF_{\frac n A}}
\to 0 
\end{align*}

\ni
As previously,  note that this conditional covariance may be written as:
\begin{align*}\label{}
& n^{1-\gamma^2}
\int_{[0,1]^2} e^{\gamma X_{\frac{n}{A}}(x)+\gamma X_{\frac{n}{A}}(y) - \gamma^2 (\ln \frac{n}{A} +1 )} \cos(2\pi n x) \sin(2\pi ny) \\
 & \hskip 3cm \times   \left ( e^{\gamma^2 K_{\frac{n}{A} } (|x-y|)} -1 \right ) dx dy\,. 
\end{align*}
The reason why real and imaginary parts are asymptotically uncorrelated is due to the trigonometric identity 
\begin{align*}\label{}
\cos(2\pi n x) \sin(2\pi ny) = \frac 1 2 \sin (4\pi n (x+y))\,,
\end{align*}
which produces cancelling  high modes in the rescaled window $|x-y|\leq A/n$ (as opposed to the  $\cos(2\pi n(x-y))$ terms that we handled for the asymptotics of $\Eb{n^{1-\gamma^2} \Re(Z_n)^2}$). 
\qed

\subsection{The case of the GMC on the unit circle.}\label{ss.circle}
$ $

The proof of Theorem \ref{thlaw} follows exactly the same lines as the proof of Theorem \ref{th.Toy}, except a more tedious white-noise type decomposition with same martingale property (i.e. from $\calF_{\frac n A}$ to $\calF$) needs to be used. 

It turns out that  the white-noise type decomposition p24 of \cite{junnila2017uniqueness} (to which we refer for details) has the analogous covariance as the covariance~\eqref{e.BB} we used above. 
Once one relies on this decomposition, the above proof applies as well to the case of the GMC on the circle.

\bibliographystyle{alpha}
\bibliography{biblio_new}

 \end{document}